\documentclass[draft,12pt,twoside,intlimits]{amsart}
\usepackage[latin2]{inputenc}
\usepackage{amsmath}
\usepackage{amsthm}
\usepackage{amssymb}
\usepackage[mathscr]{eucal}
\usepackage{enumerate}

\usepackage{amsfonts}
\usepackage{geometry,amssymb}
\newtheorem{theorem}{\bf Theorem}
\newtheorem{proposition}[theorem]{\bf Proposition}
\newtheorem{corollary}[theorem] {\bf Corollary}
\newtheorem{lemma}{\bf Lemma}

\newtheorem{definition}{\bf Definition}

\newtheorem{remark}{\bf Remark}

\newcommand{\grad}{{\rm grad\,}}

\newcommand{\con}{{\rm con~}}
\newcommand{\vol}{{\rm vol~}}

\newcommand{\NN}{\mathbb N}
\newcommand{\ZZ}{\mathbb Z}
\newcommand{\RR}{\mathbb R}
\newcommand{\CC}{\mathbb C}

\newcommand{\RRd}{{{\mathbb R}^d}}
\newcommand{\CCd}{{{\mathbb C}^d}}

\newcommand{\LL}{\mathcal L}
\newcommand{\PP}{\mathcal P}
\newcommand{\PPn}{{\mathcal P}_n}

\newcommand{\bfx}{{x}}

\newcommand{\bfy}{{y}}

\newcommand{\bfa}{{a}}

\newcommand{\bfr}{{\bf r}}

\def\sbt{\subset}

\def\la{\lambda}

\def\vp{\varphi}

\def\intt{\mathop{\rm int}}

\def\<{\langle}
\def\>{\rangle}
\def\ey{\emptyset}

\newcounter{oldresult}
\def\theoldresult{\Alph{oldresult}}

\newcounter{oldprop}
\def\theoldprop{\Alph{oldprop}}

\newcounter{oldlemma}
\def\theoldlemma{\Alph{oldlemma}}

\newcounter{oldcor}
\def\theoldcor{\Alph{oldcor}}

\newcounter{oldconjecture}
\def\theoldconjecture{\Alph{oldconjecture}}
\newenvironment{oldconjecture}{
  \em
  \vskip 0.10in
  \refstepcounter{oldconjecture}
  \noindent{\bf Conjecture\ \theoldconjecture.}
}{\vskip 0.10in}

\newcounter{hypothesis}
\def\thehypothesis{\Alph{hypothesis}}
\newenvironment{hypothesis}{
  \em
  \vskip 0.10in
  \refstepcounter{hypothesis}
  \noindent{\bf Hypothesis\ \thehypothesis.}
}{\vskip 0.10in}

\begin{document}

\title[Bernstein type estimates for multivariate polynomials]
{A comparative analysis of Bernstein type estimates for the
derivative of multivariate polynomials}

\author[Sz. Gy. R\'ev\'esz]
{Szil\'ard Gy. R\'ev\'esz}

\address{A. R\'enyi Institute of Mathematics \newline \indent Hungarian
Academy of Sciences, \newline \indent Budapest, P.O.B. 127, 1364
\newline \indent Hungary} \email{revesz@renyi.hu}


\begin{abstract} We compare the yields of two methods to obtain
Bernstein type pointwise estimates for the derivative of a
multivariate polynomial in points of some domain, where the
polynomial is assumed to have sup norm at most 1. One method, due
to Sarantopoulos, relies on inscribing ellipses into the convex
domain $K$. The other, pluripotential theoretic approach, mainly
due to Baran, works for even more general sets, and yields
estimates through the use of the pluricomplex Green function (the
Zaharjuta -Siciak extremal function). Using the inscribed ellipse
method on non-symmetric convex domains, a key role was played by
the generalized Minkowski functional $\alpha(K,x)$. With the aid
of this functional, our current knowledge is precise within a
constant $\sqrt{2}$ factor. Recently L. Milev and the author
derived the exact yield of this method in the case of the
simplex, and a number of numerical improvements were obtained
compared to the general estimates known. Here we compare the
yields of this real, geometric method and the results of the
complex, pluripotential theoretical approaches on the case of the
simplex. In conclusion we can observe a few remarkable facts,
comment on the existing conjectures, and formulate a number of
new hypothesis.
\end{abstract}



\thanks{Support of this work was rejected by Hungarian
National Foundation for Scientific Research, project \# T-047074.}
\thanks{This research was completed during the author's stay in
Sofia, Bulgaria under the exchange program of the Bulgarian and
Hungarian Academies of Sciences.}

\maketitle


\section{Introduction}\label{sec:intro}

If a univariate algebraic polynomial $p$ is given with degree at
most $n$, then by the classical Bernstein-Szeg\H o inequality
(\cite{SZ}, \cite{CS}, \cite{BO}) we have

\begin{equation}\label{classicbernstein}
|p'(x)| \leq \frac {n \sqrt{ ||p||_{C[a,b]}^2 - p^{2}({x}) }}
{\sqrt{(b-x)(x-a)}} \qquad (a<x<b).
\end{equation}

This inequality is sharp for every $n$ and every point $x\in
(a,b)$, as $$ \sup \left\{ \frac {|p'(x)|}{\sqrt{ ||p||_{C[a,b]}^2
- p^{2}({x})}} \, :\,\, \deg p \le n,\,|p(x)|< \|p\|_{C[a,b]}
\right\}= \frac {n}{\sqrt{(b-x)(x-a)}}\,\,. $$

We may say that the upper estimate \eqref{classicbernstein} is
exact, and the right hand side is just the "true Bernstein factor"
of the problem.

Polynomials and continuous polynomials are defined even on
topological vector spaces, see e.g. \cite{Din}. The set of
continuous polynomials over $X$ will be denoted by ${\mathcal P} =
{\mathcal P}(X)$ and polynomials in ${\mathcal P}$ with degree not
exceeding $n$ by ${\mathcal P}_{n} = {\mathcal P}_{n}(X)$.

In the multivariate setting a number of extensions were proved for
the classical result \eqref{classicbernstein}. However, due to the
geometric variety of possible convex sets replacing intervals of
$\RR$, our present knowledge is still not final. The exact
Bernstein inequality is known only for symmetric convex bodies,
and we are within a bound of some constant factor in the general,
nonsymmetric case.

For more precise notation we may define formally for any
topological vector space $X$, a subset $K\subset X$, and a point
$x\in K$ the $n^{\rm th}$ "Bernstein factor" as

\begin{equation}\label{bernsteinfactor}
B_n(K,\bfx):=\frac 1n \sup \left\{ \frac {\|D p(\bfx)\|}{\sqrt{
||p||_{C(K)}^2 - p^{2}({\bfx})}} \, :\,\, \deg p \le
n,\,|p(\bfx)|<||p||_{C(K)} \right\}\,,
\end{equation}
where $Dp(\bfx)$ is the derivative of $p$ at $\bfx$, and even for
an arbitrary unit vector $\bfy\in X$
\begin{equation}\label{directionalbernsteinfactor}
B_n(K,\bfx,\bfy):=\frac 1n \sup \left\{ \frac {\langle D
p(\bfx),\bfy \rangle }{\sqrt{ ||p||_{C(K)}^2 - p^{2}({\bfx})}} \,
:\,\, \deg p \le n,\,|p(\bfx)|<||p||_{C(K)} \right\}\,,
\end{equation}
where $\langle D p(\bfx),\bfy \rangle$ is the directional
derivative in direction $y$ (which equals the value attained by
the gradient, as a linear functional, on $y$).

Our aim is to investigate these and related quantities, and to
analyze the methods of estimating them.

\section{The inscribed ellipse method of Sarantopoulos}
\label{sec:inellipse}

Recall that a set $K\subset X$ is called {\it convex body} in a
normed space (or in a topological vector space) $X$ if it is a
bounded, closed convex set that has a non-empty interior. The
convex body $K$ is {\em symmetric}, iff there exists a center of
symmetry $x$ so that reflection of $K$ at $x$ leaves the set
invariant, that is, $K=-(K-x)+x=-K+2x$. In the following we will
term $K$ to be {\em centrally symmetric} if it is symmetric with
respect to the origin, i.e. if $K=-K$. This occurs iff $K$ can be
considered the unit ball with respect to a norm $\|\cdot\|_{(K)}$,
which is then equivalent to the original norm $\|\cdot\|$ of the
space $X$ in view of $B_{ X,\, \|\cdot\|}( 0,r)\sbt K\sbt
B_{X,\,\|\cdot\|}( 0,R)$.

The {\it ``maximal chord"} of $K$ in direction of $v\ne  0$ is
\begin{equation}\label{maxchord}
\begin{gathered}
\tau (K,v):= \sup\{ \la \ge 0:\exists \,\,  y,\, z\in K \,\hbox{
s.t. } z= y+\lambda v\}
\\
=\sup \big\{ \la \ge 0\colon K\cap(K+\la v)\ne \ey\big\}
\\
=\sup\{\la \ge 0\colon \la v\in K-K\}
\end{gathered}
\end{equation}
Usually $\tau(K, v)$ is not a ``maximal" chord length, but only a
supremum, however we shall use the familiar finite dimensional
terminology (see for example \cite{W}).

The {\it support function\/} to $K$, where $K$ can be an arbitrary
set, is defined for all $ v^*\in X^*$ (sometimes only for $ v^*\in
S^*:=\{ v^*\in X^*~:~ \|v^*\|=1 \} $) as
\begin{equation}\label{(2.3)}
h(K, v^*):= \sup_K v^*=\sup\big\{\< v^*, x\>\colon  x\in K\big\},
\end{equation}
and the {\it width of $K$ in direction\/} $ v^*\in X^*$ (or $
v^*\in S^*$) is
\begin{equation}\label{(2.4)}
\begin{gathered}
\hfill w(K, v^*):= h(K, v^*)+h(K,- v^*)=\sup_K v^*+\sup_K(- v^*)=
\\
=\sup\big\{\< v^*, x- y\>\colon  x, y\in K\big\}=2h\big(C,
v^*\big)=w\big(C, v^*\big).
\end{gathered}
\end{equation}
Then the {\em minimal width} of $K$ is $w(K):=\inf_{S^*}w(K,v^*)$
and the sharp inequalities
\begin{equation}\label{wtau}
w(K)\le \tau(K,v)\le {\rm diam}~ K,\quad w(K)\le w(K,v^*)\le {\rm
diam}~ K,
\end{equation}
always hold, even in infinite dimensional spaces, compare
\cite[\S 2]{RS}.

In $\RR$ the position of a point $x\in\RR$ with respect to the
"convex body" $I$ can be expressed simply by $|x|$ (as $\pm x$
occupy symmetric positions). For this in the multivariate case the
most frequent tool is the Minkowski functional. For any $ x\in X$
the {\it Minkowski functional\/} or {\it (Minkowski) distance
function\/} \cite[p.~57]{HC} or {\it gauge\/} \cite[p.~28]{Ro} or
{\it Minkowski gauge functional\/} \cite[\S 1.1(d)]{P} is defined
as
\begin{equation}\label{(1.1)}
\vp_K( x):= \inf\{\la> 0\colon  x\in\la K\}~.
\end{equation}
Clearly \eqref{(1.1)} is a norm on $ X $ if and only if the convex
body $K$ is centrally symmetric with respect to the origin. If
$K\sbt X $ is a centrally symmetric convex body, then the norm
${\|\cdot\|}_{(K)}:=\vp_K$ can be used successfully in
approximation theoretic questions as well. As said above, for
${\|\cdot\|}_{(K)}$ the unit ball of $X$ will be $K$ itself. In
case $K$ is nonsymmetric, the role of the so called {\em
generalized Minkowski functional} emerged in the above
quantification problem. This generalized Minkowski functional
$\alpha(K,x)$ also goes back to Minkowski \cite{Mi} and Radon
\cite{Rad}, see also \cite{Gr}, \cite{RS}. There are several ways
to introduce it: perhaps the shortest is the following. First let
\begin{equation}\label{gammadef}
\gamma(K, {x}) :=  \inf \left\{ 2 \frac{\sqrt {||{x} - {a}||~
||{x} - {b}||}} {||{a} - {b}||}\,: {a}, {b}
\in
\partial{K}, {\rm s.t.}~~ {x} \in [a,b] \right\}.
\end{equation}
Then we have
\begin{equation}\label{alphadef}
\alpha(K, x) = \sqrt{ 1 -\gamma^{2}(K, {x})}.
\end{equation}
In fact, usefulness of \eqref{alphadef} and the possibility of the
wide ranging applications stems from the fact that this geometric
quantity incorporates quite nicely the geometric aspects of the
configuration of $x$ with respect to $K$, which is mirrored by
about a dozen (!), sometimes strikingly different-looking,
equivalent formulations of it. For the above and many other
equivalent formulations with full proofs, further geometric
properties and some notes on the applications in approximation
theory see \cite{RS} and the references therein; for the first
appearance of it in approximation theoretical questions see
\cite{RiSh}.

The {\em method of inscribed ellipses} was introduced into the
subject by Y. Sarantopoulos \cite{Sar}. This works for arbitrary
interior points of any, possibly nonsymmetric convex body. The key
of all of the method is the next
\begin{lemma}[\bf Inscribed Ellipse Lemma, Sarantopoulos, 1991]\label{inellipse}
Let $K$ be any subset in a vector space $X$. Suppose that ${ x}
\in K$ and the ellipse
\begin{equation}\label{ellipse}
{\bf r}(t) = \cos{t}~ \bfa + b \sin{t}~ { y} + {\bfx-\bfa} \qquad
(t \in [-\pi,\pi))\,.
\end{equation}
lies inside $K$. Then we have for any polynomial $p$ of degree at
most $n$ the Bernstein type inequality
\begin{equation}\label{Bernsteinellipse}
| \langle {D} p({ x}), { y} \rangle| \leq \frac{n}{b} \sqrt{
||p||_{C(K)}^2 - p^{2}({ x})}.
\end{equation}
\end{lemma}

\begin{theorem}[\bf Sarantopoulos, 1991]\label{unball}
Let $p$ be any polynomial of degree at most $n$ over the normed
space $X$. Then we have for any unit vector $ y \in X$ the
Bernstein type inequality
\begin{equation}\label{Bernsteineq}
| \langle {D} p({ x}), { y} \rangle| \leq \frac{n \sqrt{
||p||_{C(K)}^2 - p^{2}({ x}) }} {\sqrt{1-\| x\|^{2}_{(K)}}}.
\end{equation}
\end{theorem}

\begin{theorem}[\bf Sarantopoulos, 1991]\label{scb}
Let $K$ be a symmetric convex body and $ y$ a unit vector in the
normed space $X$. Let $p$ be any polynomial of degree at most $n$.
We have $$ |\langle D p(x),y \rangle| \leq \frac{2n \sqrt{
||p_n||_{C(K)}^2 - p^{2}({ x}) }} {\tau(K,{ y}) \sqrt{1 -
\varphi^{2}(K, { x})}}~. $$ In particular, we have $$ \|D p(
x)\|\le \frac{2 n \sqrt{ ||p||_{C(K)}^2 - p^{2}({ x}) }} {w(K)
\sqrt{1-\varphi^{2}(K, { x})}}~, $$ where $w(K)$ stands for the
width of $K$.
\end{theorem}

The above solve the question for the case of a symmetric convex
body $K$. However, in the general, non-symmetric case it can be
rather difficult to determine, or even to estimate the
$b$-parameter of the "best ellipse", what can be inscribed into a
convex body $K$ through $ x\in K$ and tangential to direction of
$y$. Still, we can formalize what we want to find.
\begin{definition}[\bf Milev-R\'ev\'esz, 2003]\label{generalbestellipse}
For arbitrary $K\subset X$ and $ x\in K$, $ y\in X$ the
corresponding "best ellipse constants" are the extremal quantities
\begin{equation}\label{EKxy}
E(K, x, y):=\sup \{ b \,: \, \bfr \subset K \,\,{\hbox {\rm with}}
\,\, \bfr \,\, {\hbox {\rm as given in \eqref{ellipse}}} \} ~ .
\end{equation}
Also, for the immediately resulting estimation of the gradient, we
defined in \cite{MR}
\begin{equation}\label{EKx}
E(K, x):=\inf \{ E(K, x, y) \,: \,  y \in X , || y||=1\, \}\, .
\end{equation}
\end{definition}

Clearly, the inscribed ellipse method yields Bernstein type
estimates whenever we can derive some estimate of the ellipse
constants. In case of symmetric convex bodies, Sarantopoulos's
Theorems \ref{unball} and \ref{scb} are sharp; for the
nonsymmetric case we know only the following result.

\begin{theorem}[\bf (Kro\'o--R\'ev\'esz, \cite{KR}, 1998]\label{KrooRevesz}
Let $K$ be an arbitrary convex body, ${ x} \in {\rm int}K $ and
$\|{ y}\| = 1$, where $X$ can be an arbitrary normed space. Then
we have
\begin{equation}\label{krry}
| \langle {D}p( x),  y \rangle | \leq
 \frac {2 n \sqrt{ ||p||_{C(K)}^2 - p^{2}({ x}) }}
 { \tau(K, { y}) \sqrt{1 - \alpha(K,{ x})} } ,
 \end{equation}
for any polynomial $p$ of degree at most $n$. Moreover, we also
have
\begin{equation}\label{krrgrad}
||{D}~p({ x})|| \leq \frac {2n \sqrt{ ||p||_{C(K)}^2 - p^{2}({ x})
}} {w(K) \sqrt{ 1- \alpha(K, { x})} } \leq \frac {2 \sqrt 2 n
\sqrt{ ||p||_{C(K)}^2 - p^{2}({ x}) }} {w(K) \sqrt{ 1- \alpha^2(K,
{ x})} }.
\end{equation}
\end{theorem}

Note that in \cite{KR} the best ellipse is not found; for most
cases, the construction there gives only a good estimate, but not
an exact value of \eqref{EKxy} or \eqref{EKx}. (In fact, here we
quoted \cite{KR} in a strengthened form: the original paper
contains a somewhat weaker formulation only.)

It is worthy to recall here that geometrically the proof of
\eqref{krry} follows the following idea. To construct an ellipse
through $x$, there parallel to $y$, and inscribed into $K$ it
suffices to find the best (i.e., of maximal possible $b$
parameter) such ellipse, which is inscribed within the quadrangle
of the vertices of a {\em maximal chord in direction of $y$} (or,
in infinite dimension, some chord in the direction and
$\epsilon$-almost maximal), and the vertices of the {\em parallel
chord through the point} $x$. That ellipse is precisely
calculated, and its $b$-parameter is estimated independently of
the location of these chords (even if they degenerate into one
line, in which case the ellipse becomes a line segment). (In
general the best $b$-parameter could not be calculated, though.)
We will remind this geometrical construction later.

One of the most intriguing questions of the topic is the following
conjecture, formulated first in \cite{RS}.

\begin{oldconjecture}\label{alphasquare}{\bf
(R\'ev\'esz--Sarantopoulos, 2001).} Let $X$ be a topological
vector space, and $K$ be a convex body in $X$. For every point $ x
\in {\rm int} K$ and every (bounded) polynomial $p$ of degree at
most $n$ over $X$ we have $$ \|D p( x)\|\le \frac{2 n \sqrt{
||p||_{C(K)}^2 - p^{2}({ x}) }} {w(K) \sqrt{1-\alpha^{2}(K,
x)}}\,, $$ where $w(K)$ stands for the width of $K$.
\end{oldconjecture}

\section{Some results on the simplex}\label{sec:simplex}

We denote $|x|_2 := (\sum_{i=1}^{d} x_{i}^2)^{1/2}$ the Euclidean
norm of $x = (x_{1}, \ldots, x_{d}) \in {\mathbb R}^{d}$. Let $$
\Delta := \Delta_d := \{(x_{1},\ldots,x_{d}) : x_{i} \geq 0,
i=1,\ldots,d, \sum_{i=1}^{d}x_{i} \leq 1\} $$ be the standard
simplex in ${\mathbb R}^{d}$. For fixed $ { x} \in {\rm int}
\Delta, $ and $ { y} = (y_{1},\ldots,y_{d}),~ |{ y}|_{2} = 1$ the
best ellipse constant of $\Delta$ is, by Definition
\ref{generalbestellipse}, $ E(\Delta, x, y) $. By a tedious
calculation via the Kuhn-Tucker theorem and some geometry, the
following was obtained in \cite{MR}.
\begin{theorem}[\bf Milev-R\'ev\'esz, 2003]\label{ellipseyield}
Let $p \in {\mathcal P}_{n}^{d}$.
Then for every ${  x} \in {\rm int} \Delta$ and ${  y} \in
{\mathbb S}^{d-1}$ we have
\begin{equation}\label{directionalyield}
| D_{{  y}} p({  x}) | \leq \frac {n \sqrt{ ||p||_{C(\Delta)}^2 -
p^{2}({  x}) }} { E(\Delta,{  x}, { y}) },
\end{equation}
where $E(\Delta,x,y)$ has the precise value
\begin{equation}\label{Evalue}
E(\Delta, x, y) = \left\{ \frac{y_{1}^2}{x_{1}} + \cdots +
\frac{y_{d}^2}{x_{d}} + \frac{(y_{1} + \ldots + y_{d})^2}{1 -
x_{1} - \ldots - x_{d}} \right\}^{-1/2}.
\end{equation}
\end{theorem}

Note that the inequality
\begin{equation}\label{oldright7}
{\frac{1} {E(\Delta,{ x}, { y})} }    \leq { \frac {2}
{\tau(\Delta, { y}) \sqrt{1 - \alpha(\Delta, { x})} } }
\end{equation}
holds true for every $x \in {\rm int} \Delta$ and $y \in {\mathbb
S}^1$, which is {\em not} by chance: the general estimation of
\eqref{krry} must be valid also for $\Delta$, and the precise
value, calculated for $\Delta$, can only be better. But equality
occurs for some directions, to what we return soon.

From now on let us restrict ourselves to the case $d=2$. We denote
the vertices of $\Delta$ by $O = (0,0), A = (1,0), B = (0,1)$ and
the centroid (i.e. mass point) of $\Delta$ by $M = (1/3,1/3)$. It
is calculated in \cite{MR} that
\begin{align}\label{r}
 \alpha(\Delta, { x})  & = 1-2 r( x) \qquad \qquad \text{with} \\
  r:=r(x)& = \min \{ x_{1}, x_{2}, 1 - x_{1} - x_{2}\}
               = \left\{
           \begin{array}{ll}
           x_{1}, & {  x} \in \triangle OMB \\
            x_{2}, & {  x} \in \triangle OMA \\
           1 - x_{1} - x_{2}, & {  x} \in \triangle AMB \\
           \end{array} \right . \notag ,
\end{align}
and if ${  y} = (\cos \varphi, \sin \varphi)$ $(0 \leq \varphi
\leq \pi)$ then
\begin{equation}\label{tau}
\tau(\Delta, {  y}) = \left\{
      \begin{array}{ll}
      1/(y_{1} + y_{2}),  &  \varphi \in [0,\pi/2] \\
      1/y_{2}, & \varphi \in ( \pi/2, 3\pi/4] \\
      -1/y_{1}, & \varphi \in ( 3\pi/4, \pi]. \\
      \end{array}
      \right.
\end{equation}

Then it can be calculated that we have equality in
\eqref{oldright7} exactly for the directions $y=(\cos\vp,\sin\vp)$
with $\vp=0,\pi/2,3\pi/4$ and its $\pm$ contemporaries (i.e. at
$\vp+\pi\ZZ$ with these values of $\vp$), and for some
corresponding values of the point $x$.

Why is that so? For these, and only these vectors can we have a
coincidence of the above geometrical figure, the quadrangle
occurring in the proof of \eqref{krry}, and the exact domain into
what we must really inscribe the ellipse through $x$ and there
parallel to $y$: for all other directions the maximal chord in
direction of $y$ lies strictly inside $\Delta$ and another
ellipse, slightly stretched behind that chord, can also be
inscribed. Therefore, it is geometrically natural that it occurs
only for these directions that nothing better can be obtained,
than the ellipse calculated in Theorem \ref{KrooRevesz}, while for
other directions precise calculation of the best ellipse must
always yield a better ellipse constant.

Denote $|Dp(x)|_{2}$ the Euclidean length of the gradient vector
of $p$ at $x$ -- also equal to the operator norm $\| Dp(x)\|$
with respect to the Euclidean norm. In \cite{MR} then the
following estimates were deduced from Theorem\ref{ellipseyield}.

\begin{proposition}[\bf Milev-R\'ev\'esz, 2003]\label{gradient}
Let $p \in {\mathcal P}_{n}^{2}$. Then for every $x \in {\rm int}
\Delta$ we have
\begin{equation}\label{oldright11}
 |Dp(x)|_{2}\leq \frac{n\sqrt{||p||_{C(\Delta)}^2 - p^{2}(x)}}{E(\Delta,x)},
\end{equation}
where
\begin{equation}\label{oldright12}
E(\Delta,x) = \sqrt{\frac{2 x_{1} x_{2} (1 - x_{1} - x_{2})}
{x_{1}(1-x_{1}) + x_{2}(1-x_{2}) +  D(x)} }
\end{equation}
with
\begin{align}\label{Dxdef}
D(x):& =  \sqrt{ [x_{1}(1-x_{1}) + x_{2}(1-x_{2})]^2 - 4 x_{1}
x_{2} (1- x_{1} - x_{2})}= \\
&= \sqrt{ [x_{1}(1-x_{1}) - x_{2}(1-x_{2})]^2 + 4 x_{1}^2
x_{2}^2} > 0 \qquad \quad (\forall x \in {\rm int} \Delta)\notag
\,.
\end{align}
\end{proposition}

From this the following improvements were achieved in Theorem
\ref{KrooRevesz} for the special case of $K = \Delta$.
\begin{proposition}[\bf Milev-R\'ev\'esz, 2003]\label{alphacompare}
Let $p \in {\mathcal{P}}_{n}^{2}$ and $||p||_{C(\Delta)} = 1$.
Then for every $x \in {\rm int}~\Delta$ we have
\begin{equation}\label{oldright14}
 |{D}p(x)|_{2} \leq { \frac { \sqrt{3} n \sqrt{
||p||_{C(\Delta)}^2 - p^{2}({  x}) } }
       { w(\Delta) \sqrt{ 1 - \alpha(\Delta,x)} } }.
\end{equation}
Furthermore, using the conjectural quantity $\sqrt{ 1 -
\alpha^2(\Delta, x)}$ on the right, we even have
\begin{equation}\label{number}
 |{D}p({  x})|_{2} \leq
 \frac{ \sqrt{3+\sqrt{5}}~ n
 \sqrt{||p||_{C(\Delta)}^2 - p^{2}(x)} }
       { w(\Delta) \sqrt{ 1 - \alpha^2(\Delta, {  x})}}.
\end{equation}
\end{proposition}

The result \eqref{number} improve the constant in Theorem
\ref{KrooRevesz} but falls short of Conjecture \ref{alphasquare},
since $2\sqrt{2} = 2.8284\ldots > \sqrt{3+\sqrt{5}} =
2.2882\ldots > 2.$ On the way proving these, it was noted that no
better constants follow from the inscribed ellipse method,
interpreted so that $E(K,x)$ is considered the yield of the
ellipse method. To this we shall return.

\section{Baran's pluripotential theoretic method}\label{sec:pluripot}

Another method of considerable success in proving Bernstein (and
Markov) type inequalities is the pluripotential theoretical
approach. Classically, all that was considered only in the finite
dimensional case, but nowadays even the normed spaces setting is
cultivated. To explain these, one needs an understanding of {\em
complexifications of real normed spaces}, see e.g. \cite{MST,
BARAN4}, as well as the {\em Zaharjuta -Siciak extremal function}
$V(z)$. We start with a formulation which is perhaps easier to
digest. That is very much like the Chebyshev problem, cf. \cite[\S
8]{RS}, except that we consider it all over the complexification
$Y:=X+iX$ of $X$, take logarithms, and after normalization by the
degree, merge the information derived by all polynomials of any
degree into one clustered quantity. Namely, for any bounded
$E\subset Y$, $V_E$ vanishes on $E$, while outside $E$ we have the
definition
\begin{equation}\label{SZV}
V_E(z):=\sup \{\frac{1}{n}\log|p(z)|~:~ 0\neq p \in \PPn(Y),~
||p||_E\le 1,~n\in \NN \} \qquad (z\notin E)
\end{equation}
For $E\subset X$ one can easily restrict even to $p\in\PP(X)$.

Note that $\log|p(z)|$ is a {\em plurisubharmonic} function (PSH,
for short), as its one complex dimensional restrictions are just
logarithms of univariate polynomials over $\CC$. After
normalization by the degree, $(1/n)\log|p(z)|$ has very regular
growth towards infinity: it is at most $\log_{+}|z|+O(1)$. So it
is reasonable to consider now the Lelong class of {\em all} such
functions:
\begin{equation}\label{Lelong}
\LL (E) := \{ u~PSH ~:~u|_E \le 0,~ u(z)\le
\log|z|+O(1)~(|z|\to\infty) \}
\end{equation}
and to define
\begin{equation}\label{Lelongextremal}
U_E(z):=\sup \{ u(z)~:~ u\in\LL (E)\}~.
\end{equation}
This function may be named the pluricomplex Green function. The
Zaharjuta -Siciak Theorem says that \eqref{Lelongextremal} and
\eqref{SZV} are equal, at least as long as $E\subset \CC^d$ is
compact, which we now assume together with $E$ being a {\em
non-pluripolar} set. (A set $E\subset \CC^d$ is pluripolar, if
there exists a PSH function vanishing on $E$; if this does not
occur, the set is called non-pluripolar.) Then, as suprema of PSH
functions (subharmonic functions on all complex "lines"), they
are, modulo upper semicontinuous regularization, PSH themselves.
They play a central role in the theory.

The extension of the Laplace- and Poisson equations is the
so-called complex Monge-Ampere equation:
\begin{equation}\label{MAeq}
(\partial\overline{\partial}u)^d := d! 4^d \text{det} \left[
\frac{\partial^2 u}{\partial z_j \overline{\partial} z_k}(z)
\right] dV(z),
\end{equation}
where $dV(z)=dx_1\wedge dy_1 \wedge \dots \wedge dx_d \wedge dy_d$
is just the usual volume element in $\CCd$. At first this equation
is applied only to smooth functions $u\in PSH \cap C^2$, say, but
due to the work of Bedford and Taylor \cite{BT}, the operator
extends, in the appropriate sense, even to the whole set of
locally bounded PSH functions (which covers the case of the upper
semicontinuous regularization $V_E^{*}$ for any non-pluripolar
$E$, see e.g. \cite{Kli}). Therefore, it makes sense to consider
\begin{equation}\label{MAofV}
(\partial\overline{\partial} V_E^{*})^d ~,
\end{equation}
which is then a compactly supported measure $\lambda_E$ and is
called the {\em complex equilibrium measure} of the set $E$. It is
shown in the theory \cite{BT} that in fact the support lies in the
polynomial convex hull $\widehat{E}$ of $E$: in case $E$ is
convex, $\widehat{E}=E$ and $V_E^{*}=V_E$: moreover, it is also
shown that this measure is kind of normalized, as
$\lambda|_E(\CC^d)=\lambda|_E(\widehat{E})=(2\pi)^d$.

For the theory related to this function and some recent
developments concerning Bernstein and also Markov type
inequalities for convex bodies or even more general sets, we refer
to \cite{BARAN0, BARAN12, BARAN1, BARAN2, BARAN3, BARAN4, BCG,
BCL, BT, Kli, MAU, M, LEV, PP}.

There are further yields of the theory of PSH functions, when
applied to the Bernstein problem: here we present a few results of
Miroslav Baran. For more precise notation now we introduce
(interpreting 0/0 as 0 here)
\begin{definition}\label{gradvectorset}
\begin{equation}\label{dradset}
G(E,x):=\{\frac{\grad p(x)}{n \sqrt{\|p\|^2-p(x)^2}} ~:~ {\bf 0}
\neq p \in \PPn, n\in \NN \},
\end{equation}
and following Baran we consider also
\begin{equation}\label{conG}
\widetilde{G}(E,x):=\con G(E,x)~.
\end{equation}
\end{definition}

Clearly for any compact $E\subset \RRd$ $\sup_{n\in\NN} B_n(E,x) =
\sup_{u\in G(E,x)} \|u\|$ holds.

\begin{theorem}[\bf Baran, 1995]\label{th:equivabove}
Let $E$ be a compact subset of $\RRd$ with nonempty interior. Then
the equilibrium measure $\lambda|_E$ is absolutely continuous in
the interior of $E$ with respect to the Lebesgue measure of
$\RRd$. Denote its density function by $\lambda(x)$ for all $x\in
\intt E$. Then we have $\frac1 {d!} \lambda(x)\ge \vol
\widetilde{G}(E,x)$ for a.a. $x\in\intt E$. Moreover, if $E$ is a
symmetric convex domain of $\RRd$, then here we have exact
equality.
\end{theorem}

\begin{oldconjecture}\label{conj:Baran}{\bf (Baran, 1995).}
We have $\frac 1{d !} \lambda(x) = \vol \widetilde{G}(E,x)$ even
if $E$ is a non-symmetric convex body in $\RRd$.
\end{oldconjecture}

Now consider $E=K\subset X$, where $K$ is now a convex body. Our
more precise results in \cite{YRS}, see also \cite[\S 8]{RS},
yield $V_K(x)=\log\left(
\alpha(K,x)+\sqrt{\alpha(K,x)^2-1}\right)$. However, in the
Bernstein problem the values of $V_K$ are much more of interest
for {\em complex} points $z=x+iy$, in particular for $x\in K$ and
$y$ small and nonzero. More precisely, the important quantity is
the normal (sub)derivative
\begin{equation}\label{normalder}
D_{y}^{+} V_E(x):=\lim\inf_{\epsilon\to 0} \frac{V_E(x+i\epsilon
y)}{\epsilon}~,
\end{equation}
as this quantity occurs in the next estimation of the directional
derivative and thus also in the gradient.
\begin{theorem}[\bf Baran, 1994 \& 2004]\label{th:BaranBernstein}
Let $E\subset X$ be any bounded, closed set, $x\in \intt E$ and
$0\ne y\in X$. Then for all $p\in \PPn(X)$ we have
\begin{equation}\label{Baranestimate}
|\langle Dp(x), y \rangle|\le n D_{y}^{+} V_E(x)
\sqrt{||p||^2_E-p(x)^2}~.
\end{equation}
\end{theorem}
\begin{proof}
In fact, \cite{BARAN1} contains this only for $\RRd$ and partial
derivatives; the same estimate in case of infinite dimensional
spaces are considered in \cite{BARAN4}, but only for symmetric
convex bodies. The same estimate occurs, also without proof but
with reference to Baran, even in the recent publication
\cite{BLW}.

For a proof for arbitrary directions $y\in\RRd$ one can consider a
rotation (orthogonal transformation) $A:\RRd\to\RRd$, taking
$e_1$ to $y$, and can take $\widetilde{E}:=A^{-1}(E)$. For this
set we have $V_{\widetilde{E}}(z)=V_{A^{-1}E}(z)=V_E(A z)$.
Calculating the upper partial (sub)derivative
$\partial_1^{+}:=D_{e_1}^{+}$ we get
\begin{align}\label{rotated}
\partial_1^{+} V_{\widetilde{E}}(x) & = \lim\inf_{\epsilon\to 0+}
V_{\widetilde{E}}(x+i\epsilon e_1)= \lim\inf_{\epsilon\to 0+}
V_{E}(Ax+i\epsilon Ae_1) \notag \\ & = \lim\inf_{\epsilon\to 0+}
V_{E}(Ax+i\epsilon y) = D_{y}^{+} V_{E}(A x)~,
\end{align} so rotating back (i.e. applying the same at a rotated point
$\widetilde{x}:=A^{-1}x$) we find
\begin{equation}\label{rotbackandforth}
D_{y}^{+} V_{E}(x) = \partial_1^{+} V_{\widetilde{E}}(A^{-1}x) =
\partial_1^{+} V_{\widetilde{E}}(\widetilde{x})~.
\end{equation}
Consider now the rotated polynomial $q(x):=p(Ax)$, which is of the
same bound on $\widetilde{E}$ as $p$ on $E$, and satisfies
$Dp=D(q\circ A^{-1})= Dq\circ A^{-1} A^{-1}$ in view of the chain
rule. Hence from \eqref{rotbackandforth} and applying the estimate
\eqref{Baranestimate} only for the first partial derivative of $q$
on $\widetilde{E}$ and at an arbitrary point
$A^{-1}x=\widetilde{x} \in \intt \widetilde{E}$ (corresponding to
$x\in \intt E$) we finally obtain
\begin{align}\label{qrotated}
|\langle Dp(x), y \rangle| & =|\langle Dq(\widetilde{x}) A^{-1} ,
A e_1 \rangle| = |\langle Dq(\widetilde{x}), e_1 \rangle| \\ & \le
n D_{e_1}^{+} V_{\widetilde{E}}(\widetilde{x})
\sqrt{||q||^2_{\widetilde{E}}-q(\widetilde{x})^2} = n D_{y}^{+}
V_E(x) \sqrt{||p||^2_E-p(x)^2}~. \notag
\end{align}

\end{proof}

It is not obvious, how such theoretical estimates can be applied
to concrete cases. First, one has to find the precise value of
$V_E$, in such a precision, that even the derivative can be
computed: then the derivatives must be obtained and only then do
we really have something. However, even that is addressed by
considering the Bedford-Taylor theory of the Monge-Ampere equation
and the equilibrium measure \cite{BT}, as the density of the
equilibrium measure gives the extremal function. In some concrete
applications all that may be calculated, a particular example (see
\cite[Example 4.8]{BARAN3}) being the following.
\begin{proposition}[\bf Baran, 1995]\label{Baranexample} The extremal
function of the standard simplex in $\RRd$ is $V_{\Delta}(z)=\log
\left|h(|z_1|+|z_2|+\dots+|z_n|+|1-(z_1+z_2+\dots+z_n)|)\right|$.
Here $h(z):=z+\sqrt{z^2-1}$ is inverse to the Joukowski mapping
$\zeta \to (1/2)(\zeta+1/\zeta)$, with that choice of the
square-root that it is positive for positive $z$ exceeding 1, so
that $h$ maps to the exterior of the unit disk.
\end{proposition}

From this and the calculation with the rotated directions above,
we can calculate\footnote{The same formula is mentioned in the
recent publication \cite{BLW}, see p. 145.}
\begin{proposition} For the standard simplex $\Delta$ of $\RRd$
and with any unit directional vector $y=(y_1,\dots,y_n)$ and any
point $x=(x_1,\dots,x_n)\in \intt \Delta$ we have the formula
\begin{equation}\label{extremaldelta}
D^{+}_y V_{\Delta}(x) = \sqrt{\frac {y_1^2}{x_1}+\dots+ \frac
{y_n^2}{x_n} + \frac{(y_1+\dots y_n)^2}{(1-(x_1+\dots+x_n))}}~.
\end{equation}
\end{proposition}
\begin{proof} First we compute
$h(|z_1|+|z_2|+\dots+|z_n|+|1-(z_1+z_2+\dots+z_n)|)$ with
$z=x+i\epsilon y$. We find
\begin{align}\label{epsiloncalculus}
& h(|x_1+i\epsilon y_1|+\dots+|x_n+i\epsilon
y_n|+|1-(x_1+i\epsilon y_1+\dots+x_n+i\epsilon y_n)|) \notag \\ &=
h\left(\sqrt{x_1^2+\epsilon^2 y_1^2}+\dots+\sqrt{x_n^2+\epsilon^2
y_n^2}+\sqrt{(1-(x_1+\dots+x_n))^2+\epsilon^2 (y_1+\dots y_n)^2}
\right) \notag \\ & =
h\left(1+\epsilon^2\left[\frac{y_1^2}{2x_1}+\dots+ \frac
{y_n^2}{2x_n} + \frac{(y_1+\dots
y_n)^2}{2(1-(x_1+\dots+x_n))}\right]+O(\epsilon^4) \right) \notag
\\ &= 1 + O(\epsilon^2) + \sqrt{\left(1+ \frac{\epsilon^2}{2}
\left[\frac {y_1^2}{x_1}+\dots + + \frac {y_n^2}{2x_n} \frac{(y_1+
\dots
y_n)^2}{(1-(x_1+\dots+x_n))}\right]+O(\epsilon^4)\right)^2-1}
\notag \\ &= 1 + O(\epsilon^2) + \sqrt{\epsilon^2 \left[ \frac
{y_1^2}{x_1}+\dots+ \frac {y_n^2}{x_n} + \frac{(y_1+\dots
y_n)^2}{(1-(x_1+\dots+x_n))}\right]+O(\epsilon^4)} \notag
\\ &=
1 + \epsilon \sqrt{\frac {y_1^2}{x_1}+\dots+ \frac {y_n^2}{x_n} +
\frac{(y_1+\dots y_n)^2}{(1-(x_1+\dots+x_n))}}+ O(\epsilon^2) ~.
\end{align}
In correspondence with $|h(t)|=1$ for $t\in[-1,1]$,
$V_{\Delta}(x)=0$ for $x \in \intt \Delta$, and so we are to
calculate $D_y^{+} V_{\Delta}(x)$ by
\begin{align}\label{extremalderiv}
& D_y^{+}  V_{\Delta}(x)\notag \\ & \quad = \lim_{\epsilon\to 0+}
\frac{1}{\epsilon} \log h\left(|x_1+i\epsilon
y_1|+\dots+|x_n+i\epsilon y_n|+|1-(x_1+i\epsilon
y_1+\dots+x_n+i\epsilon y_n)| \right) \notag \\ & \quad = \sqrt{
\frac {y_1^2}{x_1}+\dots+ \frac {y_n^2}{x_n} + \frac{(y_1+\dots
y_n)^2}{(1-(x_1+\dots+x_n))}} ~ .
\end{align}
\end{proof}

Hence we are led to the following surprising corollary.
\begin{corollary}\label{cor:coincidence} The pluripotential theoretical
estimate \eqref{Baranestimate} of Baran, calculated for the
standard simplex of $\RRd$ in \eqref{extremaldelta}, gives the
result exactly identical to \eqref{directionalyield}, obtained
from the inscribed ellipse method.
\end{corollary}
Much remains to explain in this striking coincidence, the first
being the next.
\begin{hypothesis}\label{hyp:allinone} Let $K\subset X$ be a convex body.
Then for all points $x\in \intt K$ the inscribed ellipse method
and the pluripotential theoretical method of Baran results in
exactly the same estimate, i.e. for all $y\in S^*$ we have
\begin{equation}\label{ellipseeqppot}
D_{y}^{+} V_K(x) = \frac {1}{E(K,x,y)}~.
\end{equation}
\end{hypothesis}

\section{Further geometric calculations}\label{further}

At this point it seems worthy to formulate a few naturally
occurring assumptions.

\begin{hypothesis}\label{hyp:Baranright} Let $K\subset X$ be convex body.
Then for all points $x\in \intt K$ the exact Bernstein factor is
just what results from the pluripotential theoretical method of
Baran:
\begin{equation}\label{pluriexact}
B_n(K,x) = \sup_{y\in S^{+}} D_{y}^{+} V_K(x)~.
\end{equation}
\end{hypothesis}

\begin{hypothesis}\label{hyp:ellipseright} Let $K\subset X$ be convex body.
Then for all points $x\in \intt K$ the exact Bernstein factor is
just what results from  the inscribed ellipse method of
Sarantopoulos:
\begin{equation}\label{ellipseexact}
B_n(K,x) = \frac {1}{E(K,x)}~.
\end{equation}
\end{hypothesis}

These hypothesis are certainly not true for the directional
derivatives of {\em all} directions $y\in S^*$, where both methods
can be improved upon for some $y$, as is seen below. Care has to
be exercised in formulating conjectures and hypothesis in these
matters: the situation is more complex than one might like to
have, and the simple heuristics of extending the results of the
symmetric case do fail sometimes. In this respect see also
\cite{BLM, LEV, MAU} and also \cite{BLW}, where another case of
deviation from symmetric case extension is observed for the
so-called "Baran metric" on the simplex.

There is an important and immediate observation we did not utilize
until here. Namely, we found methods (actually, two equivalently
strong ones) to estimate $D_yp(x)$. However, if we are looking for
the total derivative $\grad p(x)$, then the estimate we used was
only the trivial $||\grad p(x)|| \le \sup_{y\in S^*} |D_yp(x)| $.
Can we do any better? Yes, depending on the estimating functions
we have for $D_yp(x)$, we can.

Consider e.g. the estimates from Theorem \ref{KrooRevesz}, which
was got back also for the simplex and thus the triangle $\Delta$.
For the triangle we have an explicit computation of the maximal
chords $\tau(\Delta,x)$, c.f. \eqref{tau}, and also of the
generalized Minkowski functional $\alpha(\Delta,x)$, see
\eqref{r}, so everything is explicit and we can compute the
estimating functions. As an example, consider e.g. the point
$M:=(\frac 13, \frac 13)$ and compute all quantities involved in
the normalization of the directional derivative estimates. As a
result, we can exactly determine the arising domain $H(\Delta,M)$.

In fact, it turns out that the domain $H(\Delta,M)$ what the
general estimates of Theorem \ref{KrooRevesz} describe is a
fleecy-cloud like domain which is symmetric with respect to the
origin, and its upper half is (the part above the $x$-axis of)
the union of three disks: $D \big( (\sqrt{\frac 32}, \sqrt{\frac
32}), \sqrt{3}\big) \cup D\big( (0, \sqrt{\frac 32}), \sqrt{\frac
32}\big) \cup D\big( (-\sqrt{\frac 32},0 ), \sqrt{\frac 32}
\big)$. (Here the reader may wish to draw a figure for better
visualization of the picture.) The immediate observation is that
the domain is \emph{not convex}, and so it is certainly not an
exact description of all possible directional derivatives of the
gradient.

We can conclude that if some domain
\begin{equation}\label{partialestimate}
H:=H(K,x):=\{v=t y ~:~ y=(y_1,\dots,y_d), |t| \le r(y)\}
\end{equation}
is given with $r(y)$ being a valid estimation for the directional
derivative in direction of $y$, then to find or estimate $G(K,x)$
an additional process of restricting to the "kernel" part
\begin{equation}\label{kernel}
\widetilde{H}:=\widetilde{H}(K,x):=\bigcap_{S^*} \{v ~:~ |\langle
v,\bfy\rangle|\le r(\bfy) \}
\end{equation}
is available. That is, we always have $\widetilde{G}(K,x)\subset
\widetilde{H}$. Note that $\widetilde{H}$ is a convex, symmetric
domain for whatever point set $H$.

In order to illustrate this "kernel technique", let us come back
to the  above case of estimates from Theorem \ref{KrooRevesz} for
the triangle at point $M$. After some standard considerations with
Thales circles we find that $\widetilde{H}$ is the hexagonal
domain
$$
\widetilde{H}(\Delta,M)={\rm con} \{(\sqrt{6},0),
(\sqrt{6},\sqrt{6}), (0,\sqrt{6}),(-\sqrt{6},0),
(-\sqrt{6},-\sqrt{6}), (0,-\sqrt{6}) \}~.
$$
Observe that the area of the possible stretch of $G$ is
considerably reduced from the "fleecy-cloud" domain to the derived
hexagonal domain as ${\rm area} \, H(\Delta,M) =
9+\frac{9}{2\pi}=23.137...$, while ${\rm area} \,
\widetilde{H}(\Delta, x) = 18$. For comparison recall that
Baran's Conjecture \ref{conj:Baran} would say that the area
should be $\frac 12 \lambda_{\Delta}(M) =
\frac{\pi}{\sqrt{3^{-3}}}=16.324...$.

Let us calculate the "kernel set" $\widetilde{H}(\Delta, x)$ from
the exact estimates \eqref{directionalyield},
\eqref{Baranestimate}, \eqref{extremaldelta} which we obtain from
the ellipse (and hence also of Baran's) method! We obtain the
following\footnote{These computations were executed jointly with
Nikola Naidenov from the University of Sofia during the author's
stay in Sofia in October 2004. The author regrets that in spite
of his inevitable contribution \cite{Nikola} to the work, Nikola
Naidenov chose not to be named as a coauthor.}.

\begin{proposition}\label{prop:ellipse} With the above notations,
$\widetilde{H}(\Delta, x)$ is an ellipse domain. Moreover, if the
major axis of this ellipse domain is denoted by $\mu:=\mu(x)$ and
the minor axis is denoted by $\nu:=\nu(x)$, then we have
\begin{equation}\label{ellipseaxes}
  \mu=\sqrt{\frac {2}{x_1(1-x_1) + x_2(1-x_2)+D(x)}} \quad \text{\rm and} \quad
  \nu=\sqrt{ \frac {2}{x_1(1-x_1) + x_2(1-x_2)-D(x)}},
\end{equation}
where $D(x)$ is the quantity already defined in \eqref{Dxdef}.
\end{proposition}
\begin{proof}
For any given $x\in\Delta$ our task is to solve the equation
\eqref{kernel} for the case of $K=\Delta$, the triangle, with the
estimating function $r(y)$ being the quantity
\eqref{ellipseeqppot} of $\Delta$ defined by \eqref{Evalue}. That
is, for given, fixed $x\in \Delta$ we want to determine all those
vectors $u=(u_1,u_2)\in \RR^2$, which satisfy $|\langle u, y
\rangle | \leq 1/E(\Delta,x,y)$ for all directional vectors
$y=(\cos\vp,\sin\vp)$. Writing in \eqref{Evalue} and squaring, the
defining equations present themselves as
\begin{equation}\label{Htildeeq}
\left\{ u~:~(u_1\cos \vp + u_2 \sin \vp )^2 \leq \frac {\cos^2
\vp}{x_1}+ \frac {\sin^2 \vp}{x_2}+ \frac
{(\cos\vp+\sin\vp)^2}{1-x_1-x_2} \quad (\forall \vp \in \RR)
\right\}\,.
\end{equation}
If $\cos \vp =0$, then $|\sin \vp|=1$ and the equation reduces to
\begin{equation}\label{cosnullcase}
u_2^2 \leq \frac{1}{x_2}+\frac{1}{1-x_2-x_3} =
\frac{1-x_1}{x_2(1-x_1-x_2)}~,
\end{equation}
which we need to check besides the case when $\cos^2 \vp > 0$. Put
$x_3:=1-x_1-x_2$. Division by $\cos^2 \vp > 0$ yields
\begin{equation}\label{tanfi}
(u_1 + u_2 t )^2 \leq \frac {1}{x_1}+ \frac {t^2}{x_2}+ \frac
{(1+t)^2}{x_3} \quad (\forall t:=\tan \vp \in \RR) \,,
\end{equation}
that is, writing now $z_j:=1/x_j$ for $j=1,2,3$ and ordering,
\begin{equation}\label{teq}
0\leq (z_2+z_3-u_2^2)t^2 + 2(z_3-u_1u_2) t + ( z_1 + z_3 - u_1^2)
\quad (\forall t:=\tan \vp \in \RR),
\end{equation}
which is a second degree equation in $t$. Thus the point
$u=(u_1,u_2)$ is a solution iff the discriminant of this
quadratic equation is not positive, and either the quadratic
coefficient is positive, or it is zero and then not only the
discriminant is (nonnegative hence) zero, but also the constant
term is nonnegative, too. Note that nonnegativity of the quadratic
coefficient is just \eqref{cosnullcase}, while the discriminant
condition becomes
\begin{equation}\label{discriminant}
0\leq d(x) := (z_3-u_1u_2)^2 - (z_2+z_3-u_2^2)(z_1 + z_3 -
u_1^2)\,
\end{equation}
which is again a quadratic equation, but now in the coordinates
of the point $u$. From this form a calculation leads to
$(z_2+z_3)u_1^2 + (z_1+z_3)u_2^2-2z_3u_1u_2 \leq
z_1z_2+z_1z_3+z_2z_3$, hence multiplying by $x_1x_2x_3$ and using
$x_1+x_2+x_3=1$ a few times we arrive at
\begin{equation}\label{uellipse}
a u_1^2 + b u_2^2 - c u_1 u_2 \leq 1~,
\end{equation}
where here the coefficients
\begin{equation}\label{abc}
a:=a(x):=x_1(1-x_1)\,,\quad b:=b(x):= x_2(1-x_2)\,,\quad
c:=c(x):=2x_1x_2 \,
\end{equation}
are all strictly positive. Thus \eqref{uellipse} determines an
ellipse domain indeed. Among points satisfying \eqref{uellipse} it
is not difficult to determine the maximal values of $u_2$. These
will occur at points where the function $F(u_1,u_2)$ on the left
of \eqref{uellipse} have value $1$ and a vertical gradient, i.e.
$\partial_1 F=0$, from which conditions a standard calculation
derives that $u_1=\sqrt{x_2}/\sqrt{(1-x_1)x_3}$ and
$u_2=\sqrt{1-x_1}/\sqrt{x_2 x_3}$, showing that the maximal
possible value of $u_2$ satisfies \eqref{cosnullcase} with
equality. Moreover, in this extremal case we find that also
\eqref{teq} is satisfied, the right hand side reducing to constant
$1/[x_1(1-x_1)]>0$ identically for all $t$.

It remains to determine the major and minor axes of the ellipse
domain of $u$ described by \eqref{uellipse}. In fact, the
equation is almost in a canonical form. We need only to rotate the
coordinates by
\begin{align}\label{rotation}
  v_1:&=\cos \alpha u_1 - \sin \alpha u_2 & u_1=\cos \alpha v_1 + \sin \alpha v_2
  \\
  v_2:&=\sin \alpha u_1 + \cos \alpha u_2 & u_2=\cos \alpha v_2- \sin \alpha
  v_1 \notag
\end{align}
to obtain
\begin{equation}\label{rotellipse}
A v_1^2+B v_2^2 - C v_1 v_2 \leq 1
\end{equation}
with
\begin{align}\label{ABdefined}
A:&=A(x):= a \cos^2 \alpha + b \sin^2 \alpha + c \cos \alpha \sin
\alpha\,, \notag \\
B:&=B(x):= a \sin^2 \alpha + b \cos^2 \alpha - c \cos \alpha \sin
\alpha\,, \\
C:&=C(x):= a 2 \cos \alpha \sin \alpha - b 2 \cos \alpha \sin
\alpha + c (\sin^2\alpha -\cos^2\alpha)\,\notag .
\end{align}
In case $a=b$ the rotational angle $\alpha=\pi/4$ will be proper,
as then $C$ vanishes and we get $2A=a+b+c$, $2B=a+b-c$, and
\begin{equation}\label{canellipse}
A v_1^2+B v_2^2 \leq 1.
\end{equation}
In view of the formula $D(x)=\sqrt{(a-b)^2+c^2}$, easily seen
from \eqref{Dxdef} and \eqref{abc}, \eqref{ellipseaxes} obtains
for $a=b$.

Let now $a\ne b$ and choose $\alpha= \frac 12 \arctan
\frac{c}{a-b} \in (-\pi/4,\pi/4)$, which is chosen again to
annihilate $C(x)$ and thus to reduce the equation
\eqref{rotellipse} to \eqref{canellipse}; depending on the sign
of $a-b$, the major and minor axis of the ellipse are to be
$1/\sqrt{A}$ and $1/\sqrt{B}$ or conversely. The sum of these two
axes is $a+b$, and a calculation also leads $A-B={\rm sign}\,
(a-b) \sqrt{(a-b)^2+c^2} = \pm D(x)$, hence the values
\eqref{ellipseaxes} obtain once again.
\end{proof}

So we are led to the following result.
\begin{theorem}\label{th:ellipsearea} With the above
notations, we have ${\rm area} \, \widetilde{H}(\Delta, x) =
\frac{\pi}{ \sqrt{x_1x_2(1-x_1-x_2)}}.$
\end{theorem}
\begin{proof} As is well-known, the area of an ellipse domain $E$
having major and minor axes $\mu$ and $\nu$ is ${\rm area} E = \pi
\mu \nu$, hence Proposition \ref{prop:ellipse} leads to the
asserted value.
\end{proof}

\begin{corollary}\label{Baranborder} We have $G(x)\subseteq\con G(x)\subseteq
\widetilde{H}(x)$ with ${\rm area} \, \widetilde{H}(x) = \frac12
\lambda (x)$. Hence either $\con G(x)=\widetilde{H}(x)$, for all
points $x\in\Delta$, or Baran's Conjecture \ref{conj:Baran} fails.
\end{corollary}
\begin{proof} One must compute the density function $\lambda(x)$ of
the equilibrium measure. This has already been done by Baran,
\cite[Example 4.8]{BARAN3}: we have $\lambda(x)=\frac{2\pi}{
\sqrt{x_1x_2(1-x_1-x_2)}}$\,. On comparing to Theorem
\ref{th:ellipsearea} we find the asserted identity. Since
$\widetilde{H}$ is an ellipse domain and also $ \con {G}$ is a
convex domain, $\con G(x)\subset \widetilde{H}(x)$ and equality of
their area entails that $\con G(x)=\widetilde{H}(x)$. On the
other hand if at some point $x\in\Delta$ the respective areas
differ, then ${\rm area} \,\con G(x)< {\rm area}\,
\widetilde{H}(x) = \frac12 \lambda (x)$, hence the conjectured
identity of Baran fails.
\end{proof}

\begin{remark} While using the information on the support functional
from $H(\Delta,x)$ improves upon the known area estimates, it does
not improve the maximal gradient norm estimate of \cite{MR}.
\end{remark}

Indeed, as $\widetilde{H}(\Delta,x)$ is an ellipse domain, we have
to consider the major axis of this ellipse. It turns out that in
the case of the standard triangle, this calculation yields
$\max_{v\in \widetilde{H}} \|v\| = \max_{v\in H}
\|v\|=1/E(\Delta,x)$.

Note that $\max_{v\in V} \|v\| = \max_{v\in \con V} \|v\|$ for
any set $V$, hence regarding the maximal gradient norm estimate
it makes no difference if we consider $\con G(x)$ or $G(x)$ only.
Also note that starting out from a set $H\supset G$ and
considering the "kernel" $\widetilde{H}$, we necessarily obtain a
convex set, so from $G\subset \widetilde{H}$ it follows that even
taking convex hull we still have $\con G \subset \widetilde{H}$.

\begin{corollary} Conjecture \ref{alphasquare}
and Conjecture \ref{conj:Baran} can not hold simultaneously.
\end{corollary}
\begin{proof} According to Corollary \ref{Baranborder}, Baran's
Conjecture \ref{conj:Baran} holds if only there can be no
improvement on the estimates of the ellipse (or Baran's) method
on the simplex. But than Conjecture \ref{alphasquare} fails.
Conversely, if Conjecture \ref{alphasquare} holds, then there is
an improvement at least at certain points and in certain
directions compared to the estimates of the ellipse (or Baran's)
method, hence the estimates of Corollary \ref{Baranborder} are
strictly exceeding the right value and Baran's Conjecture
\ref{conj:Baran} fails.
\end{proof}

\section{Concluding remarks}\label{s:conclusion}

Also, another real, geometric method, of obtaining Bernstein type
inequalities, due to Skalyga \cite{S1, S3}, is to be mentioned
here: the difficulty with that is that to the best of our
knowledge, no one has ever been able to compute, neither for the
seemingly least complicated case of the standard triangle of
$\RR^2$, nor in any other particular non-symmetric case the yield
of that abstract method. Hence in spite of some remarks that the
method is sharp in some sense, it is unclear how close these
estimates are to the right answer and what use of them we can
obtain in any concrete cases.

Given the above findings, it seemed to be plausible that
Conjecture \ref{alphasquare}, if not true, can be disproved by
some explicit example. To construct a polynomial with large
gradient, as compared to the norm, means to construct a highly
oscillating polynomial. For that, various natural and more
intricate ideas were tried by Nikola Naidenov \cite{Nikola} in
Sofia during the Fall of 2004. We hope he will report on his
experiences in the near future.

The author would like to thank for the enlightening comments and
suggestions of Norm Levenberg and the careful corrections of an
anonymous referee.

\end{document}